\centerline{\bf ON THERMONUCLEAR REACTION RATES}
\vskip.5cm
\noindent
\centerline{H.J. Haubold}
\centerline{UN Outer Space Office}
\centerline{Vienna International Centre, 1400 Vienna, Austria}
\vskip.3cm
\noindent
\centerline{and}
\vskip.3cm
\noindent
\centerline{A.M. Mathai}
\centerline{Department of Mathematics and Statistics}
\centerline{McGill University, 805 Sherbrooke Street West}
\centerline{Montreal, Quebec, Canada, H3A 2K6}
\vskip1cm
\noindent
\hrule
\vskip.5cm
\noindent
{\bf Abstract}
\vskip.5cm

Nuclear reactions govern major aspects of the chemical evolution of galaxies
 and
 stars. Analytic study of the reaction  rates and reaction probability
 integrals is attempted here. Exact expressions for the reaction rates and
 reaction probability integrals for nuclear reactions in the cases of
 nonresonant, modified nonresonant, screened nonresonant and resonant cases
 are given. These are expressed in terms of H-functions, G-functions and in
 computable series forms. Computational aspects are also discussed.
\vskip.3cm
\noindent
\hrule
\vskip.5cm
\noindent
{\bf 1.{\hskip.5cm}Introduction}
\vskip.5cm
Nuclear reactions govern major aspects of the chemical evolution of the 
universe
 or, at least, its building blocks: galaxies and stars. A proper understanding 
of the
 nuclear reactions that are going on in hot  cosmic plasmas, and those
 in the laboratories as well, requires a sound theory of nuclear-reaction
 dynamics (Brown and Jarmie 1990). The nuclear reaction rate is the central 
quantity connecting the
 theoretical models of universal, galactic and stellar evolution and the
 nucleosynthesizing nuclear reactions, which can be studied to some extent in
 nuclear-physics laboratories. Compilations of reaction rates and uncertainties
 (analytic expressions and tabulated temperature-dependent values) and
 astrophysical $S$-factors (analytic expressions and tabulated energy-dependent
 values) for charged-particle induced reactions are available on the
 World Wide Web  (http://csa5.lbl.gov/chu/astro.html). The rate $r_{ij}$ of
  reacting particles $i$ 
 and $j$, in the case of nonrelativistic nuclear reactions taking place in a
 nondegenerate environment, is usually expressed as
$$\eqalignno{r_{ij}&=n_in_j\left({{8}\over{\pi 
\mu}}\right)^{1\over2}\left({{1}\over{kT}}\right)^{3\over2}\int_0^{\infty}E
\sigma(E){\rm e}^{-{{E}\over{kT}}}{\rm d}E\cr
&=n_in_j<\sigma v>&(1.1)\cr}$$
where $n_i$ and $n_j$ denote the particle number densities of the reacting
 particles $i$ and $j$, $\mu={{m_im_j}\over{m_i+m_j}}$ is the  reduced mass of
 the reacting particles, $T$ is the temperature, $k$ is the Boltzmann constant,
 $\sigma(E)$ is the reaction cross section and $v=(2E/\mu)^{1\over2}$ is the
 relative velocity. Thus $<\sigma v>$ is the reaction probability integral, 
that
 is, the probability per unit time that two particles, confined to a unit
 volume, will react with each other. The basic assumptions made in (1.1) are 
that the reacting particles $i$ and $j$ have isotropic Maxwell-Boltzmann
 kinetic-energy distributions and that the kinematic of the reaction can be
 treated in the center-of-mass system, see  Fowler (1984), Mathai and Haubold
  (1988), and Imshennik (1990).
\vskip.2cm
For nonresonant nuclear reactions between nuclei of charges $z_i$ and $z_j$
 at low energies (below the Coulomb barrier), the reaction cross section has
 the form
$$\eqalignno{\sigma(E)&={{S(E)}\over{E}}{\rm e}^{-2\pi \eta(E)}\cr
\noalign{\hbox{with}}
\eta(E)&=\left({{\mu}\over2}\right)^{1\over2}{{z_iz_je^2}\over{hE^{1/2}}}\cr}$$
where $\eta(E)$ is the Sommerfeld parameter, $h$ is  Planck's quantum of
 action and $e$ is the quantum of electric charge. For a slowly varying
 cross-section factor $S(E)$ we may expand 
$$\eqalignno{S(E)&=S(0)+{{{\rm d}S(0) }\over{{\rm d}E}} E+{1\over2}{{{\rm 
d}^2S(0)}\over{{\rm d}E^2}}E^2.\cr
\noalign{\hbox{Thus}}
<\sigma v>&=\left({{8}\over{\pi \mu}}\right)^{1\over2}\sum_{\nu 
=0}^{2}{{1}\over{(kKT)^{-\nu+{1\over2}}}}{{S^{(\nu)}(0)}\over{\nu!}}\cr
&\times\int_0^{\infty}y^{\nu}{\rm e}^{-y}{\rm e}^{-zy^{-{1\over2}}}{\rm d}y\cr
\noalign{\hbox{where}}
y&={{E}\over{kT}}\hbox{ and} z=2\pi\left({{\mu}\over{2kT}}\right)^{1\over2}{{z_iz_je^2}\over{h}}.\cr
\noalign{\hbox{Then the integral to be evaluated is}}
N_{\nu}(z)&=\int_0^{\infty}y^{\nu}{\rm e}^{-y}{\rm e}^{-zy^{-{1\over2}}}{\rm 
d}y.&(1.2)\cr}$$
We will consider a general integral of the following form:
\vskip.2cm
\noindent
{\bf Nonresonant case:}
$$I_1=\int_0^{\infty}y^{\nu}{\rm e}^{-ay}{\rm e}^{-zy^{-\rho}}{\rm 
d}y,~a>0,~z>0,~\rho >0.\eqno(1.3)$$
  \vskip.2cm It may be a stringent assumption to consider a thermonuclear 
fusion plasma
 as being in thermodynamic equilibrium. If there is a cut off of the high 
energy
 tail of the Maxwell-Boltzmann distribution then  (1.2) becomes the following:
\vskip.3cm
\noindent
{\bf Nonresonant case with high energy cut-off:}
$$I_2^{d}=\int_0^{d}y^{\nu}{\rm e}^{-ay}{\rm 
e}^{-zy^{-\rho}},~d<\infty,~a>0,~z>0,~\rho >0.\eqno(1.4)$$

We consider also ad hoc modifications of the Maxwell-Boltzmann distribution for
 the evaluation of the nonresonant thermonuclear reaction rate, which acts as a
 depletion of the tail of the distribution function
 (see Kaniadakis et al. 1996) In this case the general
 integral to be evaluated is the following:
\vskip.3cm
\noindent
{\bf Nonresonant case with depleted tail:}
$$I_3=\int_0^{\infty}y^{\nu}{\rm e}^{-ay}{\rm e}^{-by^{\delta}}{\rm 
e}^{-zy^{-\rho}}{\rm d}y,~\delta >0,~z>0,~\rho >0,~a>0.~b>0.\eqno(1.5)$$

Taking into account the electron screening effects for the reacting particles
 Haubold and Mathai (1986) consider the case of the screened nonresonant
 nuclear reaction rate (see Lapenta and Quarati 1993). In this case
$$\eqalignno{r_{ij}&=n_in_j\left({{8}\over{\pi\mu}}\right)^{1\over2}\sum_{\nu=0}
^{2}{{1}\over{(kT)^{-\nu+{1\over2}}}}{{S^{(\nu)}(0)}\over{\nu!}}\cr
&\times\int_0^{\infty}y^{\nu}{\rm e}^{-y}{\rm 
e}^{-z\left(y+{{b}\over{kT}}\right)^{-{1\over2}}}{\rm d}y.&(1.6)\cr}$$
In this case we consider the following general integral:
\vskip.3cm
\noindent
{\bf Screened nonresonant case:}
$$I_4=\int_0^{\infty}y^{\nu}{\rm e}^{-ay}{\rm e}^{-z(y+t)^{-\rho}}{\rm 
d}y,~t>0,~\rho>0,~z>0,~a>0.\eqno(1.7)$$

When the nuclear cross section $\sigma(E)$ in (1.1) has a broad single 
resonance
 it can be expressed via the parametrized Breit-Wigner formula and then we 
have, see also Haubold and Mathai (1986a),
$$\eqalignno{<\sigma 
v>&=(2\pi)^{5\over2}{{z_iz_je^2R_0w\Gamma_{kl}D}\over{\mu^{1\over2}(kT)^{3\over
2}}}\cr
&\times{{1}\over{1+\left({1\over2}\Gamma_1\right)^2}}\int_0^{\infty}{{{\rm 
e}^{-ay-qy^{-{1\over2}}}}\over{(b-y)^2+g^2}}{\rm d}y&(1.8)\cr}$$
where
$$\eqalignno{a&={{1}\over{kT\left(1+\left({1\over2}\Gamma_1\right)^2\right)}},~
~b=E_{r}-{1\over4}\Gamma_0\Gamma_1,\cr
g&={1\over2}(\Gamma_0+E_{r}\Gamma_1),~~q=\bar{q}\left(1+\left({1\over2}\Gamma_1
\right)^2\right)^{1\over2},\cr
\bar{q}&=2\pi\left({{\mu}\over2}\right)^{1\over2}{{z_iz_je^2}\over{h}}=a(kT)^{1
\over2}.\cr}$$
Thus the general integral to be evaluated in this case is of the following 
form:
\vskip.3cm
\noindent
{\bf Resonant case:}
$$I_5=\int_0^{\infty}{{y^{\nu}{\rm e}^{-ay-zy^{-\rho}}}\over{(b-y)^2+g^2}}{\rm 
d}y.\eqno(1.9)$$
\vskip.2cm
In the resonant case also we can consider a modification of the
 Maxwell-Boltzmann distribution which results in a depletion of the tail. Then
 the general integral to be evaluated is the following:
\vskip.3cm
\noindent
{\bf Resonant case with the depleted tail:}
$$I_6=\int_0^{\infty}{{y^{\nu}{\rm 
e}^{-ay-by^{\delta}-zy^{-\rho}}}\over{(c-y)^2+g^2}}{\rm d}y.\eqno(1.10)$$

The aim of the present article is to give new exact analytic representations
 of the integrals in (1.2) to (1.10). Additional representations are available 
from
 Mathai and Haubold (1988). 
\vskip.5cm
\noindent
{\bf 2.{\hskip.5cm}Reduction formulae for the  reaction probability integrals}
\vskip.5cm
Writing
$$I_2^{(d)}=I_2^{(d)}(\nu,a,z,\rho)$$
we have
$$I_2^{(\infty)}=I_1.$$
Now consider
$$\eqalignno{I_3&=\int_0^{\infty}y^{\nu}{\rm e}^{-ay}{\rm e}^{-by^{\delta}}{\rm 
e}^{-zy^{-\rho}}{\rm d}y\cr
&=\sum_{m=0}^{\infty}{{(-b)^{m}}\over{m!}}\int_0^{\infty}y^{\nu+\delta m}{\rm 
e}^{-ay}{\rm e}^{-zy^{-\rho}}{\rm d}y\cr
&=\sum_{m=0}^{\infty}{{(-b)^{m}}\over{m!}}I_2^{(\infty)}(\nu+\delta 
m,a,z,\rho)&(2.1)\cr}$$
where $(a)_n$ denotes the Pochhammer symbol,
$$(a)_n=a(a+1)...(a+n-1),~(a)_0=1,~a\ne 0.$$
Note that
$$\eqalignno{I_4&=\int_0^{\infty}y^{\nu}{\rm e}^{-ay}{\rm 
e}^{-z(y+t)^{-\rho}}{\rm d}y\cr
&={\rm e}^{at}\int_{t}^{\infty}(u-t)^{\nu}{\rm e}^{-au}{\rm 
e}^{-zu^{-\rho}}{\rm d}u\cr
&={\rm 
e}^{at}\sum_{m=0}^{\infty}{{(-\nu)_m}\over{m!}}t^m\int_{u=t}^{\infty}u^{-m}{\rm 
e}^{-au}{\rm e}^{-zu^{-\rho}}{\rm d}u\cr
&={\rm 
e}^{at}\sum_{m=0}^{\infty}{{(-\nu)_m}\over{m!}}t^m\left[I_2^{(\infty)}(-m,a,z,\
rho)-I_2^{(t)}(-m,a,z,\rho)\right].&(2.2)\cr}$$
For simplifying $I_5$ and $I_6$ we will use the identity
$${{1}\over{(c-y)^2+g^2}}=\int_0^{\infty}{\rm e}^{-[(c-y)^2+g^2]x}{\rm d}x$$
and rewrite the single integral as a double integral. That is,
$$\eqalignno{I_5&=\int_0^{\infty}{{y^{\nu}{\rm e}^{-ay}{\rm 
e}^{-zy^{-\rho}}}\over{(c-y)^2+g^2}}{\rm d}y\cr
&=\int_{x=0}^{\infty}{\rm e}^{-g^2x}\int_{y=0}^{\infty}y^{\nu}{\rm 
e}^{-x(c-y)^2}{\rm e}^{-ay}{\rm e}^{-zy^{-\rho}}{\rm d}y{\rm d}x.\cr}$$
Now we expand
$$\eqalignno{{\rm 
e}^{-x(c-y)^2}&=\sum_{m=0}^{\infty}{{(-1)^m(c-y)^{2m}x^m}\over{m!}}\cr
&=\sum_{m=0}^{\infty}\sum_{m_1=0}^{2m}{{2m}\choose 
m_1}(-1)^{m+m_1}{{c^{2m-m_1}}\over{m!}}x^{m}y^{m_1}.\cr}$$
The integral over $x$ gives
$$\int_{x=0}^{\infty}x^m{\rm e}^{-g^2x}{\rm d}x=(g^2)^{-(m+1)}m!.\eqno(2.3)$$
Substituting  back we have
$$\eqalignno{I_5&={{1}\over{g^2}}\sum_{m=0}^{\infty}\sum_{m_1=0}^{2m}{{(2m)!c^{
2m-m_1}}\over{m_1!(2m-m_1)!}}{{(-1)^{m+m_1}}\over{(g^2)^m }}\cr
&\times\int_0^{\infty}y^{\nu+m_1}{\rm e}^{-ay}{\rm e}^{-zy^{-\rho}}{\rm d}y\cr
&={{1}\over{g^2}}\sum_{m=0}^{\infty}\sum_{m_1=0}^{2m}{{2m}\choose 
m_1}{{(-1)^{m_1}}\over{c^{m_1}}}\left(-{{c^2}\over{g^2}}\right)^m I_2^{ 
(\infty)}(\nu+m_1,a,z,\rho).&(2.4)\cr}$$
Thus it is seen that all the integrals $I_1$ to $I_6$ can be reduced to the
 integral $I_2^{(d)}(\nu,a,z,\rho)$ for two different situations of 
non-negative
 as well as negative $\nu$. We will evaluate $I_2^{(d)}$ in the next section
 by using a statistical technique.
\vskip.5cm
\noindent
{\bf 3.{\hskip.5cm}Evaluation of the  integral $I_2^{(d)}$}
\vskip.5cm
In general, integrals $I_1$ to $I_6$ are quite difficult to evaluate
 analytically. Here we
 will use a statistical technique. We will evaluate the density of a product of
 two independently distributed real scalar random variables by using two
 different methods, one procedure leading to the integral that we want to
 evaluate and the other procedure leading to a representation in terms of a
 known function. Then appealing to the uniqueness of the density we have the
 integral evaluated in terms of a known special function. Let $x$ and $y$ be 
real scalar
 random variables having the densities
$$f_1(x)=\cases{c_1{\rm e}^{ -ax},~0<x<d\cr
0,~{\hbox{elsewhere}}\cr}$$
and
$$f_2(y)=\cases{c_2y^{\nu}{\rm e}^{-zy^{\rho}},~0<y<\infty\cr
0,~{\hbox{elsewhere}}\cr}$$
where $c_1$ and $c_2$ are normalizing factors such that
$$\int_{x=0}^{d}f_1(x){\rm d}x=1\hbox{  and  }\int_{y=0}^{\infty}f_2(y){\rm 
d}y=1.$$
Since the variables are assumed to be independently distributed the joint
 density of $x$ and $y$ is the product of $f_1(x)$ and $f_2(y)$. Let us
 transform $x$ and $y$ to $u=xy$ and $v=x$. Then the joint density of $u$ and
 $v$, denoted by $g(u,v)$, and the marginal density of $u$, denoted by 
$g_1(u)$,
 are given by
$$\eqalignno{g(u,v)&=c_1c_2u^{\nu}v^{-\nu-1}{\rm e}^{-av}{\rm 
e}^{-cv^{-\rho}},~c=zu^{\rho}\cr
\noalign{\hbox{and}}
g_1(u)&=c_1c_2u^{\nu}\int_0^{d}v^{-\nu-1}{\rm e}^{-av}{\rm e}^{-cv^{-\rho}}{\rm 
d}v,~ c=zu^{\rho}.\cr}$$
Hence we have
$$\int_0^{d}v^{-\nu-1}{\rm e}^{-av}{\rm e}^{-cv^{-\rho}}{\rm 
d}v={{u^{-\nu}}\over{c_1c_2}}g_1(u),~c=zu^{\rho}.\eqno(3.1)$$
Let us look at the $(s-1)$-th moment of $u$. Due to independence
$$E(u^{s-1})=[E(x^{s-1})][E(y^{s-1})].$$
But
$$\eqalignno{E(x^{s-1})&=c_1\int_0^{d}x^{s-1}{\rm e}^{-ax}{\rm d}x\cr
&=c_1\sum_{m=0}^{\infty}{{(-1)^m(ad)^m}\over{m!}}{{d^s}\over{s+m}},~\hbox{ for 
}~d<\infty\cr
&=c_1a^{-s}\Gamma(s), \Re (s)>0,~~\hbox{ for }~d=\infty\cr
\noalign{\hbox{where $\Re (\cdot)$ denotes the real par of $(\cdot)$, and}}
E(y^{s-1})&=c_2\int_0^{\infty}y^{\nu+s-1}{\rm e}^{-zy^{\rho}}{\rm d}y\cr
&={{c_2}\over{\rho}}z^{-{{\nu+s}\over{\rho}}}\Gamma\left({{\nu 
+s}\over{\rho}}\right),~~\hbox{  for  }\Re (\nu+s)>0.\cr}$$
Taking the inverse Mellin transform of $E(u^{s-1})$ we have
$$\eqalignno{g_1(u)&={{c_1c_2}\over{\rho}}z^{-{{\nu}\over{\rho}}}\sum_{m=0}^{\i
nfty}{{(-1)^m(ad)^m}\over{m!}}\cr
&\times {{1}\over{2\pi 
i}}\int_{L}{{d^s}\over{s+m}}z^{-{{s}\over{\rho}}}\Gamma\left({{\nu+s}\over{\rho
}}\right)u^{-s}{\rm d}s&(3.2)\cr}$$
where $L$ is a suitable contour and $i=\sqrt{-1}$. This contour integral
 can be written as an H-function, see for example Mathai and Saxena (1978).
 That is,
$${{1}\over{2\pi 
i}}\int_{L}{{d^{s}}\over{s+m}}z^{-{{s}\over{\rho}}}\Gamma\left({{\nu 
+s}\over{\rho}}\right)u^{-s}{\rm 
d}s=H_{1,2}^{2,0}\left[{{uz^{{{1}\over{\rho}}}}\over{d}}\bigg\vert^{(m+1,1)}_
{(m,1),\left({{\nu}\over{\rho}},{{1}\over{\rho}}\right)}\right].\eqno(3.3)$$
Substituting (3.3) in (3.2) and then comparing with (3.1) we have
$$\eqalignno{\int_0^{d}v^{-\nu-1}{\rm e}^{-av}{\rm e}^{-zu^{\rho}v^{-\rho}}{\rm 
d}v&=I_2^{(d)}(-\nu-1,a,zu^{\rho},\rho)&(3.4)\cr
&={{z^{-{{\nu}\over{\rho}}}u^{-\nu}}\over{\rho}}\sum_{m=0}^{\infty}{{(-ad)^m}\o
ver{m!}}\cr
&\times 
H_{1,2}^{2,0}\left[{{uz^{{{1}\over{\rho}}}}\over{d}}\bigg\vert^{(m+1,1)}_{(m,1)
,\left({{\nu}\over{\rho}},{{1}\over{\rho}}\right)}\right],~\hbox{ for  }
d<\infty&(3.5)\cr
&={{z^{-{{\nu}\over{\rho}}}u^{-\nu}}\over{\rho}}H_{0,2}^{2,0}\left[uaz^{{1}\over
{\rho}}\bigg\vert_{(0,1),\left({{\nu}\over{\rho}},{{1}\over{\rho}}\right)}
\right],~\hbox{  for }d=\infty&(3.6)\cr}$$
where $\Re (\nu)>0, \Re (a)>0, \Re (z) >0, \Re (\rho)>0, 0\le u\le\infty$. 
Note that
(3.4) , (3.5) and (3.6) give an $I_2^{(d)}(\nu,a,z,\rho)$ with $\nu$ negative.
 If the integral for a positive $\nu$ is required then we proceed as follows:
Take $f_1(x)=c_3x^{\nu}{\rm e}^{-ax}$ and $f_2(y)=c_4{\rm e}^{-zx^{\rho}}$,
 where $c_3$ and $c_4$ are the new normalizing constants, and then proceed as
 before. Then we end up with the following representations.
$$\eqalignno{\int_0^{d}v^{\nu -1}{\rm e}^{-av}{\rm e}^{-zv^{-\rho}}{\rm 
d}v&=I_2^{(d)}(\nu-1,a,z,\rho)&(3.7)\cr
&={{d^{\nu}}\over{\rho}}\sum_{m=0}^{\infty}{{(-ad)^m}\over{m!}}\cr
&\times 
H_{1,2}^{2,0}\left[{{z^{{{1}\over{\rho}}}}\over{d}}\bigg\vert\matrix{(\nu 
+m+1,1)& \cr
(\nu +m,1),&\left(0,{{1}\over{\rho}}\right)\cr}\right]\hbox{for} d<\infty&(3.8)\cr
&={{a^{-\nu}}\over{\rho}}H_{0,2}^{2,0}\left[az^{{{1}\over{\rho}}}|_{(\nu,1),
\left(0,{{1}\over{\rho}}\right)}\right],\hbox{  for  } d=\infty&(3.9)\cr}$$
where $\Re (\nu)>0,~\Re (a)>0,~\Re (z)>0, ~\Re (\rho)>0$. Thus (3.7), (3.8)
 and (3.9) cover the case of a positive $\nu$ in $I_2^{(d)}(\cdot)$.
 When $\rho$ is real and rational then the H-functions appearing in (3.3) to
 (3.9) can be reduced to  G-functions by using the multiplication formule for 
gamma functions, namely,
$$\Gamma(mz)=(2\pi)^{{{(1-m)}\over2}}m^{mz-{1\over2}}\Gamma(z)\Gamma\left(z+{{1}
\over{m}}\right)...\Gamma\left(z+{{m-1}\over{m}}\right),~m=1,2,....\eqno(3.10)$$
For the theory and applications of G-functions see for example Mathai (1993).
 For $\rho ={{m}\over{n}},~m,n=1,2,...$ reduction to the G-function is 
available
 from Mathai and Haubold (1988). 
\vskip.2cm
The parameters of interest in nuclear astrophysics are 
$I_2^{(d)}(\nu,a,z,\rho)$ for
 $a=1,~z>0,~\rho ={1\over2}$. In this case computable representations will be
 discussed in the next section.
\vskip.5cm
\noindent
{\bf 4.{\hskip.5cm}Computable series representations of the reaction rate
 integrals}
\vskip.5cm
Let us start with (3.9) for $\rho ={1\over2}$. Then the H-function to be
 evaluated is the following:
$$\eqalignno{H_{0,2}^{2,0}(\cdot)&={{1}\over{2\pi i}}\int\Gamma(\nu 
+s)\Gamma(2s)(az^2)^{-s}{\rm d}s\cr
&={{1}\over{2\sqrt{\pi}}}{{1}\over{2\pi 
i}}\int\Gamma(s)\Gamma\left(s+{1\over2}\right)\Gamma(\nu 
+s)\left({{az^2}\over4}\right)^{-s}{\rm d}s&(4.1)\cr}$$
by expanding $\Gamma (2s)$ with the help of (3.10). Note that (4.1) is a
 G-function of the type $G_{0,p}^{p,0}(\cdot)$, see
 for example Mathai (1993).
 Observe that for $\nu \ne \pm$
 $ {{\lambda}\over2},$
 $~ \lambda =0,1,...$
all the poles
 of the integrand are simple. Then (4.1) generates a series corresponding to
 each gamma in the integrand. Corresponding to $\Gamma (s)$ the poles are at
 $s =-n,~n=0,1,...$ and the corresponding residue is
$$\eqalignno{\lim_{s\rightarrow -n}\left[(s+n)\Gamma(s)\Gamma\left(s+{1\over2}
\right)\Gamma(\nu+s)\left({{az^2}\over4}\right)^{-s}\right]&={{(-1)^n}\over{n!}}
\Gamma\left(-n+{1\over2}\right)\Gamma(\nu -n)\left({{az^2}\over{4}}\right)^n.\cr
\noalign{\hbox{But}}
\Gamma\left(-n+{1\over2}\right)&={{(-1)^n\Gamma\left({1\over2}\right)}\over
{\left({1\over2}\right)_n}}\cr
\noalign{\hbox{and}}
\Gamma(\nu -n)&={{(-1)^n\Gamma(\nu)}\over{(1-\nu)_n}}.\cr}$$
Hence the sum of the residues gives
$$\Gamma\left({1\over2}\right)\Gamma(\nu)\sum_{n=0}^{\infty}{{(-1)^n\left({
{az^2}\over4}\right)^n}\over{\left({1\over2}\right)_n(1-\nu)_n}}=\Gamma\left
({1\over2}\right)\Gamma(\nu){_0F_2}\left(~;{1\over2},1-\nu;-{{az^2}\over4}\right)\eqno
(4.2)$$
where ${_0F_2}(\cdot)$ is a hypergeometric series which is convergent for all
 $a$ and $b$. Thus (3.9) is a linear function of 3 such series of ${_0F_2}$'s
 for $\nu\ne\pm{{n}\over2},~n=0,1,...$. If $\nu$ is an integer or half-integer
 then we have one set of poles of order 2 each and one set of order one each.
 Techniques for handling the integral when the integrand has higher order poles
 are given in Mathai (1993).  Series representations of (4.1) for all cases of
 $\nu$ are given in Mathai and Haubold (1988).
\vskip.2cm
Now let us examine (3.8) for $\rho ={1\over2}$. In this case the H-function to
 be evaluated is given by
$$\eqalignno{H&=H_{1,2}^{2,0}\left[{{z^2}\over{d}}\bigg\vert^{(\nu +m+1,1)}_{
(\nu+m,1),(0,2)}\right]\cr
&={{1}\over{2\pi i}}\int 
{{1}\over{\nu+m+s}}\Gamma(2s)\left({{z^2}\over{d}}\right)^{-s}{\rm d}s\cr
&={{1}\over{2\pi i}}\int {{1}\over{\nu 
+m+s}}\Gamma(s)\Gamma\left(s+{1\over2}\right)\left({{z^2}\over{d}}\right)^{-s}{
\rm d}s.&(4.3)\cr}$$
At $s=-\nu-m$ there is a pole of order one if 
$\nu\ne\pm{{\lambda}\over2},~\lambda =0,1,...$,
 otherwise it will be a pole of order 2 at this point. When it is a pole of
 order one the residue is given by
$$\eqalignno{{{1}\over{2\sqrt{\pi}}}&\Gamma(-\nu -m)\Gamma\left(-\nu 
-m+{1\over2}\right)\left({{z^2}\over{d}}\right)^{m+\nu}\cr
&={{1}\over{2\sqrt{\pi}}}{{\Gamma(-\nu)\Gamma\left(-\nu-{1\over2}\right)}\over{
(\nu 
+1)_m\left(\nu+{1\over2}\right)_m}}\left({{z^2}\over{d}}\right)^{m+\nu}.\cr}$$
Hence corresponding to this residue the term in (3.8) is the following:
$$\eqalignno{{{z^{2\nu}}\over{\sqrt{\pi}}}&\Gamma(-\nu)\Gamma\left(-\nu+{1\over
2}\right)\sum_{m=0}^{\infty}{{1}\over{(\nu+1)_m\left(\nu+{1\over2}\right)_m}}(-
az^2)^m\cr
&={{z^{2\nu}}\over{\sqrt{\pi}}}\Gamma(-\nu)\Gamma\left(-\nu+{1\over2}\right)
{_0F_2}\left(~;\nu+1,\nu+{1\over2};-az^2\right)\cr}$$
which is evidently convergent.
\vskip.2cm
At $s=-n,~n=0,1,...$ the integrand in (4.3) has poles of order one when
 $\nu\ne\pm{{\lambda}\over2},~\lambda =0,1,...$ The sum of the residues here
 is given by
$$\eqalignno{{{1}\over{2\sqrt{\pi}}}&\sum_{n=0}^{\infty}{{1}\over{\nu+m-n}}{{(-
1)^n}\over{n!}}\Gamma\left(-n+{1\over2}\right)\left({{z^2}\over{d}}\right)^n\cr
&={1\over2}\sum_{n=0}^{\infty}{{1}\over{\nu+m-n}}{{1}\over{\left({1\over2}
\right)_n n!}}\left({{z^2}\over{d}}\right)^n.&(4.4)\cr}$$
The term corresponding to this in (3.8) is the following double series:
$$d^{\nu}\sum_{m=0}^{\infty}\sum_{n=0}^{\infty}{{(-ad)^m}\over{m!}}{{1}\over
{\nu+m-n}}{{1}\over{\left({1\over2}\right)_n}}{{(z^2/d)^n}\over{n!}}.$$
By Horn's theorem on convergence, see for example Srivastava and Karlsson
 (1985, pp. 56--57) this double series is evidently convergent for all $a,z$ 
and
 $d$. The third term of (3.8) as well as all terms for other cases of $\nu$ can
 be easily seen to give convergent series for all $0<d<\infty,~a>0,~z>0$.
 Similar arguments hold good for (3.5) and (3.6) also.
\vskip.2cm
Let us examine $I_3$. We have expressed $I_3$ in terms of $I_2^{(\infty)}$ in 
(2.1). That is,
$$I_3=\sum_{m=0}^{\infty}{{(-b)^m}\over{m!}}I_2^{(\infty)}(\nu +\delta 
m,a,z,\rho).
$$For $\rho ={1\over2}$ a typical term in this $I_2^{(\infty)}$ behaves like a 
${_0F_2}$ given in (4.2). A typical term will be a constant multiple of a 
double series of the form
$$\sum_{m=0}^{\infty}\sum_{n=0}^{\infty}{{(-b)^m}\over{m!}}{{1}\over{\left({1
\over2}\right)_n(1-\nu-\delta 
m)_n}}{{\left(-{{az^2}\over4}\right)^n}\over{n!}}.\eqno(4.4)$$
This by Horn's theorem is convergent for all $(x,y)|0<x<\infty,~0<y<\infty$
 where $x=b$ and $y={{az^2}\over4}$. Hence the series representation of $I_3$ 
is
 convergent for all $a,b,z$ and $\delta$ for $\nu\ne\pm{{\lambda}\over2},$
$~\lambda=0,1,...$. When $\nu$ is an integer or half-integer then the series
 representation will contain psi functions and logarithmic terms but from the
 structure in (4.4) one can see that the series will be convergent.
\vskip.2cm
Now let us examine the complicated form coming from $I_5$ or from (2.4). 
Here we have the factor $I_2^{(\infty)}(\nu +m_1,a,z,\rho)$ with $m_1\ge 0$.
 For the case $\rho ={1\over2}$ we have from (3.9) and (4.1)
$$\eqalignno{\int_0^{\infty}x^{\nu+m_1}{\rm e}^{-ax}{\rm e}^{-zx^{-{1\over2}}}{\rm 
d}x&=I_2^{(\infty)}\left(\nu+m_1,a,z,{1\over2}\right)\cr
&={{a^{-(\nu+m_1+1)}}\over{\sqrt{\pi}}}{{1}\over{2\pi 
i}}\int_{L}\Gamma(s)\Gamma\left(s+{1\over2}\right)\Gamma(\nu+m_1+1+s)\left({{az
^2}\over4}\right)^{-s}{\rm d}s&(4.5)\cr}$$
Writing (4.5) as a G-function, see Mathai (1993), we have
$$I_2^{(\infty)}\left(\nu +m_1,a,b,{1\over2}\right)={{a^{-(\nu 
+m_1+1)}}\over{\sqrt{\pi}}}
 G^{3,0}_{0,3}\left[{{az^2}\over4}\bigg\vert_{0,{1\over2},\nu 
+m_1+1}\right].\eqno(4.6)$$
The behavior of this G-function for small and large values of ${{az^2}\over4}$
 is available from Mathai and Saxena (1973, p. 307) or Luke (1969, pp. 
178--180). For small values of
 ${{az^2}\over4}$ this G-function behaves like unity and hence
 $I_2^{(\infty)}$ behaves like ${{a^{-(\nu+m_1)}}\over{\sqrt{\pi}}}$. For large
 values of ${{az^2}\over4}$ the $I_2^{(\infty)}$ in (4.6) behaves like$$
\left({{az^2}\over4}\right)^{{1\over6}}\left({{az^2}\over{4}}\right)^{{{\nu 
+m_1}\over3}}{\rm e}^{-3\left({{az^2}\over4}\right)^{1\over3}}.$$
Hence for checking the convergence of $I_5$ in (2.4) it is sufficient to 
examine 
the two series 
$$\eqalignno{\sum_{m=0}^{\infty}\sum_{m_1=0}^{2m}&{{2m}\choose 
m_1}\left(-{{1}\over{c}}\right)^{m_1}a^{-m_1}\left(-{{c^2}\over{g^2}}\right)^m\
cr
&=\sum_{m=0}^{\infty}\left[\left(1-{{1}\over{ca}}\right)^2\right]^m\left(-{{c^2}
\over{g^2}}\right)^m\cr
&=\left[1+{{c^2}\over{g^2}}\left(1-{{1}\over{ca}}\right)^2\right]^{-1}~~\hbox{for}
~~\left\vert{{c}\over{g}}\left(1-{{1}\over{ca}}\right)\right\vert <1\cr}$$ 
or for  ${{1}\over{a}}-|g|<c<{{1}\over{a}}+|g|$ which is equivalent to 
 $0<c<1+|g|$ for $|g|\ge 1,$ when $a=1$,  and
$$\eqalignno{\sum_{m=0}^{\infty}\sum_{m_1=0}^{2m}&{{2m}\choose 
m_1}\left(-{{1}\over{ac}}\right)^{m_1}\left[\left({{az^2}\over{4}}\right)^{1
\over3}\right]^{m_1}\left(-{{c^2}\over{g^2}}\right)^m\cr
&=\sum_{m=0}^{\infty}\eta^m\left(-{{c^2}\over{g^2}}\right)^m=\left[1+\eta{{c^2}
\over{g^2}}\right]^{-1},\cr}$$
where 
$$\eta =\left[1-{{1}\over{c}}\left({{z}\over{2a}}\right)^{2\over3}\right]^2$$
for $\left\vert\eta{{c^2}\over{g^2}}\right\vert <1$ or for
$$2a[c-|g|]^{3\over2}<z<2a[c+|g|]^{3\over2}.$$
For example, combined with the condition for small values, we have
 $0<z<2(2+|g|)$ for $a=1,~|g|\ge 2$.

\vskip.2cm
For other values of $\rho$ also the procedure remains the same. We may 
 observe that instead of the $I_5$ in (1.9) we can also consider a more general
 integral of the type
$$I_7=\int_0^{\infty}{{y^{\nu}{\rm 
e}^{-ay-zy^{-\rho}}}\over{\left[(c-y)^2+g^2\right]^{d}}}{\rm 
d}y,~\nu,~a,~z,~c,~\rho,~d~>0.\eqno(4.7)$$
In this case replace
$${{1}\over{\left[(c-y)^2+g^2\right]^{d}}}={{1}\over{\Gamma(d)}}\int_0^{\infty}
x^{d -1}{\rm e}^{-x[(c-y)^2+g^2]}{\rm d}x,~\Re (d)>0.\eqno(4.8)$$
Then the modification in (2.3) becomes
$$\int_0^{\infty}x^{m+d -1}{\rm e}^{-g^2x}{\rm 
d}x=(g^2)^{-(m+d)}\Gamma(m+d)\eqno(4.9)$$
and a corresponding expression for (2.4) is obtained. Convergence conditions
 can be checked exactly the same way as in the case of $d -1=0$.
\vskip.2cm
$I_6$ can be reduced to a form corresponding to (2.4). Convergence conditions
 can be checked by converting the series form into one of the standard triple
 series discussed in Srivastava and Karlsson (1985) or by using the general
 procedure for a multiple series or by writing the kernel function as a
 G-function, as described above, and then checking the behavior for large and
 small values of the argument of this G-function. Note that $I_7$ can be
 expressed in terms of $I_2^{(\infty)}$ as follows:
$$\eqalignno{I_7&=\sum_{m=0}^{\infty}(-1)^m{{(d)_m}\over{m!}}(g^2)^{-(m+d)}\sum
_{m_1=0}^{2m}{{2m}\choose m_1}(-1)^{m_1}c^{2m-m_1}\cr
&\times\sum_{n=0}^{\infty}{{(-b)^n}\over{n!}}I_2^{(\infty)}(\nu+m_1+\delta 
n,a,z,\rho).&(4.10)\cr}$$
Writing $I_2^{(\infty)}$ in terms of a G-function and then checking the
 behavior of the G-function for small and large values of the argument one can
 verify the existence of $I_7$. For small values, the G-function in $I_2$ 
behaves
 like unity and then for $\rho ={1\over2}$
$$\eqalignno{I_7&=\sum_{m=0}^{\infty}(-1)^m{{(d)_m}\over{m!}}(g^2)^{-(m+d)}\sum
_{m_1=0}^{2m}{{2m}\choose m_1}(-1)^{m_1}c^{2m-m_1}\cr
&\times\sum_{n=0}^{\infty}{{(-b)^n}\over{n!}}{{a^{-(\nu+1+m_1+\delta 
n)}}\over{\sqrt{\pi}}}\cr
&={{a^{-(\nu+1)}}\over{\sqrt{\pi}}}{\rm 
e}^{-a^{-\delta}}g^{-2d}{_1F_0}\left(d;~;-{{c^2}\over{g^2}}\left(1-{{1}\over{ca}}
\right)^2\right)\cr}$$
for ${{c^2}\over{g^2}}\left[1-{{1}\over{ca}}\right]^2<1$ or
 ${{1}\over{a}}-|g|<c<{{1}\over{a}}+|g|$. For large values the behavior of the 
G-function is available from (5.2) later on and in this case the conditions for
 the existence are given in (5.8) later on.
\vskip.2cm
When $\rho ={{m}\over{n}},~m,n=1,2,...$ it is easy to see that the
 H-function in all the integrals $I_1$ to $I_6$ reduce to G-functions, of 
course
 with more parameters. General series representations of all forms of 
G-functions are available from Mathai (1993).
\vskip.5cm
\noindent
{\bf 5.{\hskip.5cm}Computational aspects}
\vskip.5cm
Anderson, Haubold and Mathai (1994) looked into the computational aspects of 
the
 integrals $I_1$ to $I_4$. Exact computations and graphs are given there for
 $\rho={1\over2},~a=1$ in the cases of $I_1$ for $\nu =1,2$; $I_2$ for
 $\nu =0,~d=1,~5$; $I_3$ for $(\nu,\delta,b)=(0,2,0.001),~ (1,5,1)$; $I_4$
 for $(\nu,t)=(0,1),~ (0,5)$. For large values of $z$ one can use the
 asymptotic forms of the integrals for computational purposes. These can be
 worked out by using the asymptotic form of the G-function. From (3.9) and 
(4.1)
$$I_2^{(\infty)}\left(\nu,a,z,{1\over2}\right)={{a^{-(\nu+1)}}\over{\sqrt{\pi}}}
G_{0,3}^{3,0}\left[{{az^2}\over4}\bigg\vert_{\nu+1,a,{1\over2}}\right].\eqno(5.1)$$
But  for large values of ${{az^2}\over4}$ we have 
$$G_{0,3}^{3,0}\left[{{az^2}\over4}\bigg\vert_{\nu 
+1,0,{1\over2}}\right]\approx {{2\sqrt{\pi}}\over{\sqrt{3}}}a^{-(\nu +1)}{\rm 
e}^{-3\left({{az^2}\over4}\right)^{1\over3}}\left({{az^2}\over4}\right)^{{{2\nu 
+1}\over6}}.\eqno(5.2)
$$By substituting this expression for the G-function in $I_1$ to $I_7$ we get
 the following forms for large values of ${{az^2}\over4}$ with $a=1,~\rho 
={1\over2}$:$$
\eqalignno{I_1&\approx 
2\left({{\pi}\over3}\right)^{1\over2}\left({{z^2}\over4}\right)^{{{2\nu 
+1}\over6}}{\rm e}^{-3\left({{z^2}\over4}\right)^{1\over3}}&(5.3)\cr
I_2&\approx d^{\nu +1}{\rm 
e}^{-d}\left({{z^2}\over{4d}}\right)^{-{1\over2}}{\rm 
e}^{-2\left({{z^2}\over{4d}}\right)^{1\over2}}&(5.4)\cr
I_3&\approx 
2\left({{\pi}\over3}\right)^{1\over2}\left({{z^2}\over4}\right)^{{{2\nu+1}\over
6}}{\rm e}^{-3\left({{z^2}\over4}\right)^{1\over3}}{\rm 
e}^{-b\left({{z^2}\over4}\right)^{{{\delta}\over3}}}&(5.5)\cr
I_4&\approx 2\left({{\pi}\over3}\right)^{1\over2}{\rm 
e}^{t}\left({{z^2}\over4}\right)^{1\over6}{\rm 
e}^{-3\left({{z^2}\over4}\right)^{1\over3}}\left[\left({{z^2}\over4}\right)^{1\
over3}-t\right]^{\nu}&(5.6)\cr
I_5&\approx{{2}\over{g^2}}\left({{\pi}\over3}\right)^{1\over2}{\rm 
e}^{-3\left({{z^2}\over4}\right)^{1\over3}}\left({{z^2}\over4}\right)^{{{2\nu 
+1}\over6}}\left[1+\eta{{c^2}\over{g^2}}\right]^{-1}&(5.7)\cr
\noalign{\hbox{where $\eta 
=\left[1-{{1}\over{c}}\left({{z}\over2}\right)^{2\over3}\right]^2$}}
I_6&=I_7\hbox{  for  }d=1\cr
I_7&\approx{{2\sqrt{\pi}}\over{\sqrt{3}}}a^{-(\nu+1)}g^{-2d}\left({{az^2}\over4}
\right)^{{1\over2}+{{\nu}\over3}}{\rm 
e}^{-3\left({{az^2}\over4}\right)^{1\over3}}\cr
&\times {\rm 
e}^{-{{b}\over{a^{\delta}}}\left({{az^2}\over4}\right)^{{{\delta}\over3}}}{_1F_0}
(d;~;-\gamma^2),\cr}$$
for $|\gamma|<1$ where 
$$\gamma^2={{c^2}\over{g^2}}\left[1-{{1}\over{c}}\left({{z}\over{2a}}\right)^{2
\over3}\right]^2.$$
For a broad overview of the aspects of numerical evaluation of special
 functions, including available software packages for this purpose, see on the
 World Wide Web: http://math.nist.gov/nest/.
\vskip.5cm
\noindent
\centerline{\bf References}
\vskip.5cm
\noindent
Anderson, W.J., Haubold, H.J. and Mathai, A.M. (1994). Astrophysical
 thermonuclear functions. {\it Astrophysics and Space Science}, {\bf 214},
 49--70.
\vskip.5cm
\noindent
Brown, R.E. and Jarmie, N. (1990). Differential cross sections at low energies
 for $^2H(d,p)^3H$ and $^2H(d,n)^3He$. {\it Phys. Rev.}, {\bf C41}, 1391--1400.
\vskip.5cm
\noindent
Fowler, W.A. (1984). Experimental and theoretical nuclear astrophysics:
 the quest for the origin of the elements. {\it Rev. Modern Phys.}, {\bf 56},
 149--179.
\vskip.5cm
\noindent
Haubold, H.J. and Mathai, A.M. (1986). Analytic results for screened 
non-resonant nuclear reaction rates. {\it Astrophysics and Space Science},
 {\bf 127}, 45--53.
\vskip.5cm
\noindent
Haubold, H.J. and Mathai, A.M. (1986a). Analytic representations of
 thermonuclear reaction rates. {\it Studies in Applied Mathematics.},
 {\bf LXXV (2)}, 123--137.
\vskip.5cm
\noindent
Imshennik, V.S. (1990). Nuclear reaction rates for interacting particles with
 anisotropic velocity distributions (in Russian). {\it Sov. Plasma Phys.},
 {\bf 16}, 655--663
\vskip.5cm
\noindent
Kaniadakis, G., Lavagno, A. and Quarati, P. (1996). Generalized statistics and
 solar neutrinos. {\it Phys. Lett.}, {\bf B369}, 308--312.
\vskip.5cm
\noindent
Lapenta, G. and Quarati, P. (1993). Analysis of non-Maxwellian fusion reaction
 rates with electron screening. {\it Z. Phys.}, {\bf A346},  243--250.
\vskip.5cm
\noindent
Luke, Y.L. (1969). {\it The Special Functions and Their Approximations, Vol.I},
 Academic Press, New York.
\vskip.5cm
\noindent
Mathai, A.M. (1993). {\it A Handbook of Generalized Special Functions for
 Statistical and Physical Sciences}, Oxford University Press, Oxford.
\vskip.5cm
\noindent
Mathai, A.M. and Haubold, H.J. (1988). {\it Modern Problems in Nuclear and
 Neutrino Astrophysics}, Akademie-Verlag, Berlin.
\vskip.5cm
\noindent
Mathai, A.M. and Saxena, R.K. (1973). {\it Generalized Hypergeometric Functions
 with Applications in Statistics and Physical Sciences}, Springer-Verlag, 
(Lecture Notes in Mathematics No.348),
 Heidelberg.
\vskip.5cm
\noindent
Mathai, A.M. and Saxena, R.K.(1978). {\it The H-function with Applications in
 Statistics and Other Disciplines}, Wiley New York.
\vskip.5cm
\noindent
 Srivastava, H.M. and Karlsson, P.W. (1985). {\it Multiple Gaussian
 Hypergeometric Series}, Ellis-Horwood, Chichester, U.K.

\bye